\documentclass[11pt, twoside, a4paper]{article}
\usepackage{amsfonts, amssymb, amsthm}
\usepackage[latin1]{inputenc}
\usepackage[intlimits,fleqn]{amsmath}
\usepackage[english]{babel}
\usepackage{graphicx}

\include{epsf}

\newcommand{\la}{\langle}
\newcommand{\ra}{\rangle}

\newcommand{\suli}[2]{\sum\limits_{#1}^{#2}}

\def\be{\begin{eqnarray*}}
\def\bel{\begin{eqnarray}}
\def\ee{\end{eqnarray*}}
\def\eel{\end{eqnarray}}
\def\a{\alpha}

\def\w{\omega}
\def\wb{\bar{\omega}}
\def\bw{\bar{\omega}}
\def\wh{\widehat}

\def\N{\mathbb{N}}
\def\E{\mathbb{E}}

\def\Cn{C_{\bar{o}}}

\def\bE{\bar{E}}

\def\bare{\bar{e}}

\def\Z{\mathbb{Z}}

\def\T{\mathcal{T}}

\def\eb{\bar{e}}

\def\h{{\bf h}}

\def\sme{\setminus^e}

\def\eqiv{\Leftrightarrow}

\def\vec{\left(\begin{array}{cc} }
\def\cev{\end{array} \right)}

\theoremstyle{plain}
\newtheorem{theo}{Theorem}[section]
\newtheorem{lemma}[theo]{Lemma}

\theoremstyle{definition}

\begin{document}
\begin{titlepage}

{\center
{\large 

\quad  {\bf Amenability of horocyclic products of percolation trees}} 

\vspace{0.5cm}

Florian Sobieczky \footnote{Institut f\"ur Mathematik C, TU-Graz, supported by
  FWF (Austrian Science Fund), project P18703}\\ 
\vspace{0.5cm}}

\vspace{2cm}

{\bf Abstract:} For horocyclic products of percolation subtrees of
regular trees, we show almost sure amenability. Under a symmetry 
condition concerning the growth of the two percolation trees,
we show the existence of an {\em increasing} F\o lner sequence
(which we call {\em strong} amenability).\\



\end{titlepage}  
\newpage

\section{Introduction}

\subsection{Amenability and strong amenability of Diestel-Leader graphs}

The family of graphs called {\em Diestel-Leader graphs} (=: DL graphs), which
are of exponential growth, have the special property that they include an
amenable subfamily, and are otherwise non-amenable and even non-unimodular (see
[\ref{cakawo}]). They have been designed to answer a question posed by Woess
[\ref{sowo}] about the existence of quasi-transitive graphs which are not
quasi-isometric to a Cayley graph of a group. An answer has been given recently
in [\ref{efw}]. These graphs are certain (`horocyclic') products (for the
definition, see section \ref{sec:horo}) of homogeneous trees, which, if taken
for two trees of equal degree, turn out to be the Cayley graphs of the
lamplighter group on the integers [\ref{barwoe}, \ref{cheper}]. As shown by
Kaimanovich and Vershik [\ref{kv}], the speed of the simple random walk on the
Cayley graph of the lamplighter group on $\Z^d$ is zero, iff $d\in\{1,\;2\}$.
Therefore, the speed of simple random walk on DL-graphs for trees of equal
degree is zero.\\ 

Virag[\ref{vir}] has proven that the positivity of the speed of a simple random
walk on an infinite graph is implied by the positivity of the anchored
isoperimetric constant ({\bf anchored expansion}).  It can be shown, that for
symmetric Diestel-Leader graphs (i.e. trees of equal degree), there is a F\o
lner sequence [\ref{barwoe}], implying amenability. Even the anchored
isoperimetric constant vanishes, a situation which is called {\bf strong
  amenability} in [\ref{hss}]. From a paper by Chen, Peres and Pete this
follows for Bernoulli percolation on symmetric Diestel-Leader graphs
[\ref{cheper}].  Our results, however, refer to a modification of Bernoulli
percolation: similar to Diestel-Leader graphs, two trees are `coupled' by
requiring their Busemann functions [\ref{woess}] to add up to zero. However,
the tees involed are independent bond-percolation subtrees of regular trees.
In spite of there possibly being no vertices of degree two, leading to large
subgraphs consisting only of finite chains (`stretchings'), there is
amenability, almost surely. Moreover, when the process is chosen in a symmetric
way on the factors of the horocyclic product, we show a.s. {\em strong}
amenability (Theorem \ref{theo:strong_amen}).\\

It has been shown in [\ref{cheper}] (see also [\ref{lyo}]) that the super-critical
percolation cluster of an invariant percolation containing a
pre-assigned `root' is a.s.  weakly non-amenable in the case of an
underlying non-amenable graph and the percolation being Bernoulli if
the retention parameter $p$ is sufficiently close to one.  In
[\ref{hss}], it has been proven that this cluster remains a.s.
strongly amenable under the `random pertubation' given by invariant
percolation, if the underlying graph is amenable and transitive.  Our
results complement these findings by showing amenability of a certain
percolation model different from Bernoulli percolation.  This
percolation process is defined by coupling two independent subtrees
which result from independent bond-percolation. The coupling is done
according to the rules of constructing Diestel-Leader graphs
([\ref{woess2}], see also [\ref{woess}], chapter 12.18), which we will
call the horocyclic product.  We prove a.s.  strong amenability at a
point $p_o$ in the range of parameters at which there is equal growth
of the random trees involved.

\subsection{Notation}

The letter $G=\la V, E\ra$ will be denoted fore the deterministic
`underlying' transitive (Diestel-Leader) graph, on the edges of which
a percolation process will be defined. `$H(\bw)$' (or just $H$) will
be reserved for the random (`percolative') subgraphs $\la \bar{V},
\bar{E}\ra$ of $G$. Since we will deal with a product probability
space $\bar{\Omega}=\Omega'\times\Omega$, its elements will be called
$\bw=\la \w', \w\ra$, throughout. Edges $\eb \in \bar{E}$ will be
undirected and denoted by subsets of the vertices: $\bar{e}=\{\bar{k},
\bar{l}\}\subset V$.  There are no loops, such that every $\bar{e}\in
E$ has two elements. The set $\Cn\subset \bar{V}$ will be the
connected component of $H$ containing a pre-assigned root $\bar{o}\in
V$. We will focus on bond percolation graphs, such that $V=\bar{V}$
and $H$ is a so called partial graph of $G$.

\subsection{Products of trees with a fixed end}

We recall the definitions concerning {\em trees with a fixed end}.  For the
following definitions, we refer to [\ref{cakawo}] and [\ref{woess}], for a more
detailed discussion.  A ray is an infinite sequence of successive neighboured
vertices without repetitions. In a homogeneous tree $\T_M$ of degree $M=q+1\ge
3$, denote by $\partial \T_M$ its {\em boundary}, which is the union of all {\em
  ends}. An end is an equivalence class of {\em rays}, where two rays are
equivalent if both have infinitely many vertices in common with a third. In
particular, for trees, this means that the traces of two rays of the same end
differ only by finitely many vertices.\\

After having chosen a {\em root}, denoted by `$o$', and an element $\gamma$ of
$\partial\T_{q+1}$, it is possible to define the {\em Busemann function}
$\h(x):= d(x,c_\gamma)-d(c_\gamma,o)$, where $x\curlywedge \gamma$ is the {\em
  confluent}, the last common vertex of the two geodesic rays between $\gamma$
and $x$ and $\gamma$ and $o$, and $d(x,y)$ is the length of the geodesic ray
between $x$ and $y$.  $\h(x)$ is the index of the `level' of the vertex $x$ on
the directed tree with fixed end $\gamma$ (see [\ref{vadim}] , where the term
pointed tree is used).  In Diestel-Leader graphs, this level hierarchy is used
to construct a product of trees (see [\ref{woess}], chapter 12.18).  If this
product involves two homogeneous trees of equal degree as its factors, the
automorphism group of the resulting graph is the Cayley graph of an amenable
group [\ref{woess2}]. This results from the group being a closed subgroup of
the Cartesian product of the amenable automorphism groups of the two involved
trees. On the other hand, if the trees are homogeneous of {\em different}
degree, the graph is not even unimodular ([\ref{woess}], chap. 12.18).

\subsection{Horocyclic products}\label{sec:horo}

Let $\T'=\la V(\T'), E(\T') \ra$, and $\T=\la V(\T), E(\T) \ra$ be two
homogeneous trees of degree $\alpha'+1$ and $\alpha+1$ with fixed roots $o'$
and $o$, and fixed ends $\gamma'$ and $\gamma$, respectively.  Due to the fixed
end and fixed root, the vertex $k\in V(\T)$ has a Busemann-function
(=:level-coordinate) $\h(k)$, likewise for $k'\in V(\T')$: $\h'(k')$. Let
$DL_{\a',\a}= \la V, E\ra$ be the graph with

\bel
V &=& \{ \;\la k', k \ra \in V(\T')\times V(\T)\;\;\;|\;\;\; \h(k') = -
\h(k)\;\},\label{eq:V}  
\eel

while the edge-set $E$  is inherited by the edge-sets $E(\T')$ and
$E(\T)$:
\bel
 E = \{ \; \{\la k', k \ra , \; \la l', l \ra \}\subset V\;|\; k'\sim
l', k\sim l \},
\eel

with $\sim$ meaning neighbours in $\T'$ and $\T$. The graphs $DL_{\a',\a}$ with
$\a, \a' \in \{1,2, 3, 4, ... \}$ are the Diestel Leader graphs. Since it is a
subgraph of the product of $\T'$ and $\T$, we denote this {\em horocyclic
  product} by
\bel
\T' \circ \T := G := DL_{\a', \a}.
\eel

Let $V_h\subset V$ be the vertex set of the finite, connected induced subgraph
$G^{(h)}$ of $G$, which contains all vertices $\bar{v}:= \la v', v\ra$ with

\be
|\h(v')| \le h,\;\;\;\textrm{ and }\;\;\;\;|\h(v)| \le h,
\ee

for some given $h\in\N$ (see Fig. 1).\\

To conclude the introductory section, we present a lemma about horocyclic
products, which will be used in the proofs of the main results. For two graphs
$G_1$ and $G_2$ with disjoint vertex sets, we call $G_1 + G_2 := \la V(G_1)
\cup V(G_2), E(G_1) \cup E(G_2) \ra$ the graph consisting of the connected
components given by those of $G_1$ and $G_2$ without there being any additional
connection (edge) between any of them ({\em disjoint union}, see [\ref{wes}]).
By $G_1 \cup_{E'} G_2$ we mean the graph $\la V(G_1) \cup V(G_2), E(G_1)\cup
E(G_2) \cup E'\ra$, such that if each of $G_1$ and $G_2$ is connected and
$e=\{k, l\}$ with $k\in V(G_1), l\in V(G_2)$, then $G_1 \cup_{\{e\}} G_2$ is
connected.  We say {\bf $E'$ connects $G_1$ and $G_2$}.

\begin{lemma}\label{lemma:first} Let $T', T$ be two infinite trees, for each of which a Busemann
  function has been defined. Let $T_1, T_2$ subtrees of $T$ with disjoint
  vertex sets, then $T'\circ (T_1 + T_2)=T'\circ T_1 + T'\circ T_2$, and this
  graph is disconnected. If $E'\subset E(T)$ connects $T_1$ and $T_2$, then
  $T'\circ (T_1 \cup_{E'} T_2)$ is connected.
\end{lemma}

{\em Proof:} Let $k\in V(T_1)$, $l\in V(T_2)$, then $k,l \in V(T_1 + T_2)$ with
no path in $T_1 + T_2$ connecting $k$ and $l$. Since any path in $T'\circ (T_1
+ T_2)$ connecting $\la k', k\ra$ with $\la l', l\ra$ for some vertices $k', l'
\in V(T_1)$ is a graph of the form $P'\circ P$ with $P'$ a path in $T'$
(connecting $k'$ and $l'$, and $P$ a path in $T$ connecting $k$ and $l$,
$T'\circ (T_1 + T_2)$ is not connected, since such a path $P$ does not exist.
On the other hand, given two vertices $\la k', k_1\ra$ and $\la l', k_2\ra $
with $k_1, k_2 \in V(T_1)$, two paths can be found, $P'\le T'$, connecting $k'$
and $l'$, and $P\le T_1$, connecting $k_1$, and $k_2$, since $T'$ and $T_1$ are
connected. Then the graph $P'\circ P$ is a connected subgraph of $T'\circ T_1$,
and since all edges are distinct, a connecting path itself, namely, connecting
$\la k', k_1\ra$ and $\la l', k_2\ra $. The same holds for $T'\circ T_2$.
Since every vertex in $T'\circ (T_1 + T_2)$ is either in $V(T' \circ T_1)$ or
$V(T'\circ T_2)$, the graph $T'\circ (T_1 + T_2)$ consists of exactly these two
connected components.  Finally, when $E'$ connects $T_1$ and $T_2$, for every
two vertices $k, l \in V(T_1 \cup_{E} T_2)$, there is a path $P\le T_1 \cup_{E}
T_2$ connecting them, such that any two vertices in $V(T'\circ(T_1 \cup_{E}
T_2))$ can be connected by a path of the form $P'\circ P$, where $P'$ is a
connecting path of $k'$ and $l'$ in $T'$. \hfill \qed

\section{Strong amenability of random horocyclic 
products}\label{sec:strong}

{\bf Definition:} For any finite subgraph $H_f=\la V_f, E_f\ra$ of $G=\la V,
E\ra$ of order $|V_f|$, define the {\bf isoperimetric ratio $I_G(H_f)$ of $H_f$ in $G$} by

\bel
I_G(H_f)\;:=\; \frac{|\partial_G V_f|}{|V_f|}\label{eq:iso_r}
\eel

where $|\cdot|$ is cardinality, and $\partial_G V_f=\{ k\in V\setminus V_f\;
:\;\{k,l\}\in E, l\in V_f\;\}$ the (outer vertex-) boundary of $V_f$ in $G$.
The {\bf anchored isoperimetric constant} [\ref{tho}][\ref{vir}] (in [\ref{vir}],
$|\cdot|$ denotes the volume= sum of weights of a subgraph) is given by

\bel 
I:= \liminf\limits_n\{\;I_G(H)\;|\; H \textrm{ is an order-$n$
connected subgraph of $G$ containing vertex }\bar{o}\;\}.
\label{eq:def_I}
\eel

Note, that a slightly different definition of $I$ can be given, which
depends on the choice of the root (compare with the appendix of
[\ref{hss}]).  However, positivity of either constant implies the
positivity of the other.  We say, that $G$ is {\bf strongly amenable},
if $I$ is zero.  Otherwise, $G$ is {\bf weakly non-amenable}, or, that
it has {\bf anchored expansion}[\ref{vir}].\\

Recalling the definition of an amenable graph, which is given if the
 {\bf isoperimetric constant} $I_o(G)=\inf\{ I_G(H)| H\le G, H $
finite$\}$ is zero, it is clear that an infinite graph with a vanishing {\em
  anchored} isoperimetric constant
is amenable.\\

\begin{figure}[htb]
\centerline{
\epsfxsize=10cm
\epsffile{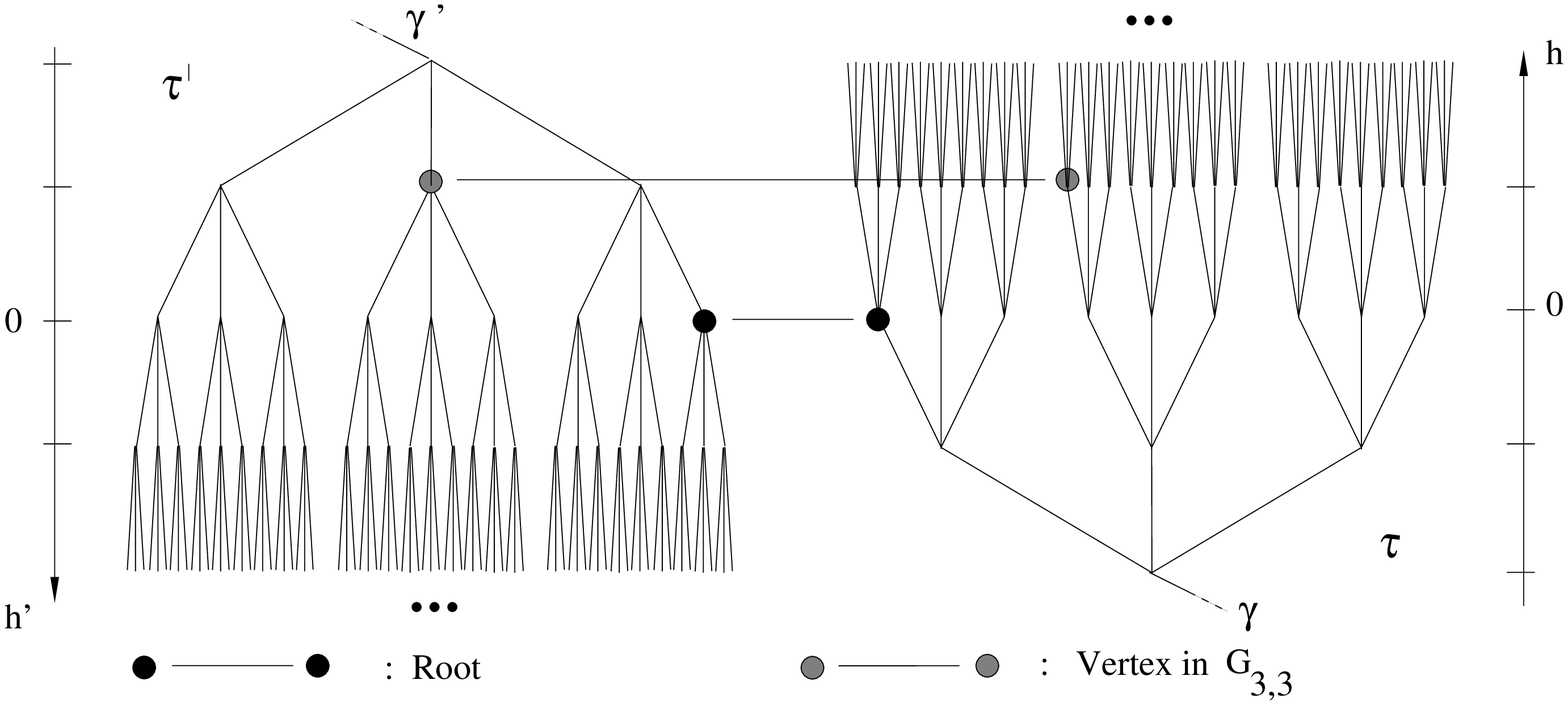} \hskip 1.5cm}
\caption{Subgraph $G^{(h)}$ of $G=\la V, E\ra=DL_{\a', \a}$ with $h=2, \a=\a'=3$
(The horizontal bars and the upside-down representation of one of the two trees
refer to the condition of definition (\ref{eq:V}), [\ref{woess2}]. An edge are
two horizontal bars with vertices which are neighbours in the trees.)} 
\end{figure}

Now, we consider a specific bond-percolation $\mu$ on the product sigma-algebra
of $\Omega = 2^{E(\T)}$ and likewise a percolation $\mu'$ on the product
sigma-algebra of $ \Omega' = 2^{E(\T')}$. Among the edges of $G$, choose a `set
of {\bf remanent edges} $E_r\subset E$ and call $E_p:= E \setminus E_r$ the `set of
{\bf percolative edges}', in the following way.  For a realization $\bw =\la \w',
\w\ra \in\Omega'\times\Omega$ in the product-probability space, let $H(\bw)=\la
V, \bE(\bw)\ra$ be the partial graph of $G$, given by

\be 
\bE(\bw)&:=& \; \{\; \;\{ \la k', k \ra, \la l', l\ra \} \in
E_p\;\;|\;\;\w'(\{k', l'\})=1,\;\w(\{k, l\})=1\;\} \;\cup \;E_r.   
\ee

$H(\bw)$, for $\bw\in\Omega'\times\Omega$, is a bond-percolative
subgraph of the Diestel-Leader graph $G$. At first, to increase
clarity, we will be interested in the concrete example of $\a=\a'=3$,
and in the special Bernoulli percolation which allows only taking away
`the third edges'. For this purpose, for all vertices $\bar{k}=\la k',
k \ra\in V$, we {\bf mark} exactly one of the three edges of $\T$ (and
$\T'$) incident to the vertices $k$ (and $k'$) pointing away from the
fixed end $\gamma$ (and $\gamma'$) to be the set of pairs of edges
(percolative edges $E_p$) which contain at least one marked edge of either
tree. Let $E_r$ (remanent edges) be the complement of $E_p$ with
respect to $E$.\\

\begin{figure}[htb]
\centerline{
\epsfxsize=9cm
\epsffile{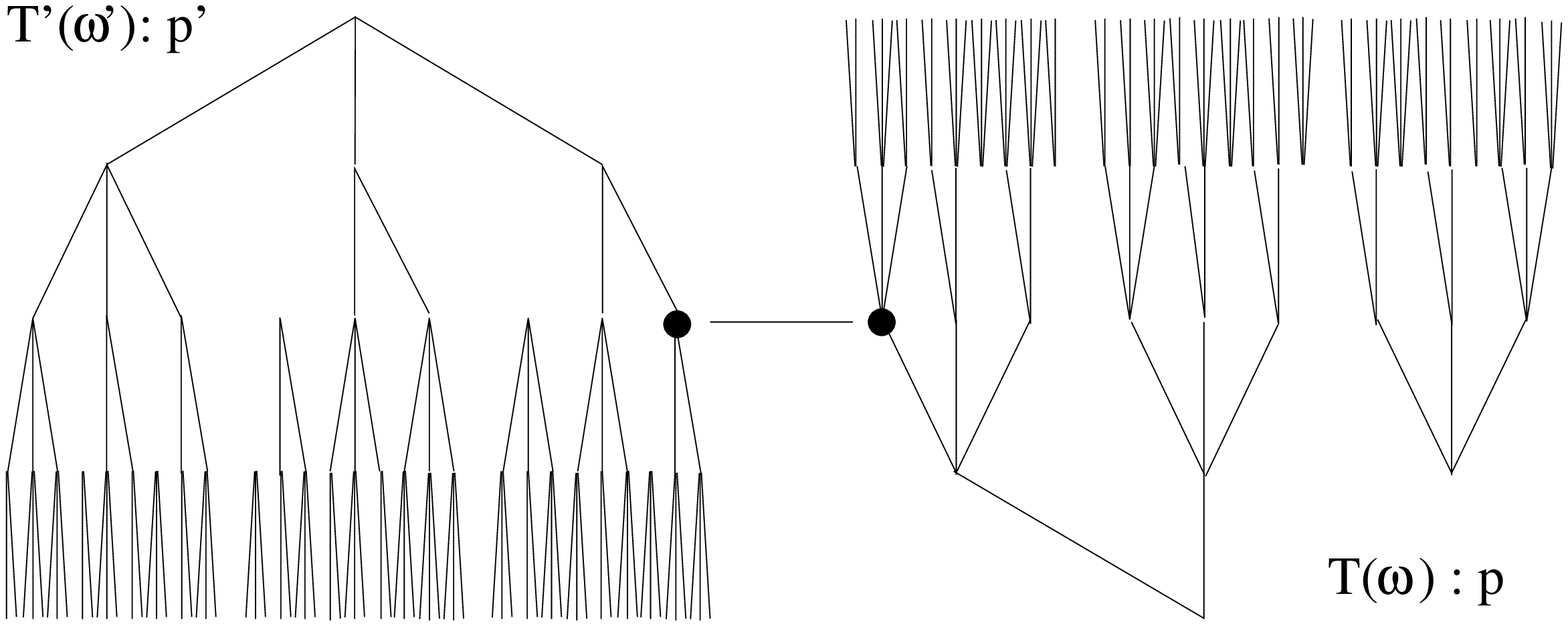} \hskip 1.5cm}
\caption{Percolative subgraph $H^{(2)}(\bw)$ of $G^{(2)}$:
  $\alpha_o'=\alpha_o=2$, and $\a'=\a=3$}
\end{figure}

In other words, for an edge $\bar{e}=\{ \la k', k\ra, \la l', l\ra\}\in E$ to
be a remanent edge ($\in E_r$), neither of the two edges $\{k', l'\}\in E(\T')$
and $\{k, l\}\in E(\T)$ may be marked. The edges in $E_r$ are not subject to
the removal of the percolation. Let $H^{(h)}(\bw)$ be the subgraph of $H(\bw)$
with vertices $\bar{v}=\la v', v\ra$ fulfilling $|\h'(v')|\le h, |\h(v)|\le h$,
then Figure 2 shows a typical realization of $H(\bw)$, when the unmarked edges
consist of the first two children, and  `only some of the third
  edges are removed.'\\

  We define the percolation measure $\bar{\mu}_{\bare}$, with $\bare\in E$ by
  the following: Let $E_m(\T)$ (respectively, $E_m(\T')$) be the {\bf marked
    edges} of $\T$ (respectively, $\T'$), and $E_u(\T)$ (resp. $E_u(\T')$) be
  the {\bf unmarked} ones: choose $E_u(\T)$ by enumerating all edges incident to
  any vertex and pointing away from the end $\gamma$. Then call each of them
  part of $E_u(\T)$, if the number labeling it is smaller or equal to a fixed
  integer $\alpha_o\in\{1, ..., \alpha\}$. Do the same to define $E_u(\T')$,
  with the integer $\alpha_o'$ as the number of `unmarked children'. Let 

 \bel
  \mu_e(\w) =\chi_{E_u(\T)}(e)\;+\;\chi_{E_m(\T)}(e)
  \;\left(\;p\cdot\chi_{\{\w(e)=1\}}(\w)\;+\;(1-p)\cdot
    \chi_{\{\w(e)=0\}}(\w)\right), \label{def:mu}
\eel

with some $p\in(0,1)$, and $\mu'_{e'}$, with $e'\in E(T'), p'\in(0,1)$ by 

\bel
\mu_{e'}(\w') =\chi_{E_u(\T')}(e')\;+\;\chi_{E_m(\T')}(e')
\;\left(\;p'\cdot\chi_{\{\w'(e')=1\}}(\w')\;+\;(1-p')\cdot
  \chi_{\{\w'(e')=0\}}(\w')\right). \label{def:mup}
\eel

Let $\bar{\mu}: \bare = \la e', e\ra \in E(\T')\times E(\T) \mapsto
\mu_{e'}\otimes \mu_{e}$ be the product measure on the product
sigma-algebra of $\Omega'\times\Omega$. Figure 2 shows a subgraph of a
realization with $\alpha_o'=\alpha_o=2, \;\alpha=\alpha'=3$. Note that
we have for $\bar{e}=e'\circ e$ that $\bar{\mu}(\bar{e} \textrm{
  open})=\mu'_{e'}(e' \textrm { open})\mu_e(e \textrm{
  open})$. However, $\bar{\mu}$ is not an independent percolation
measure, since e.g. the events that the two different edges $e'\circ
e_1$ and $e'\circ e_2$ be open are positively correlated.\\

Then, for $\bw \in \Omega'\times\Omega$, and $p',p\in [0,1]$, call $\Cn(\bw)$
the connected component of $H(\bw)$, containing a preassigned root $\bar{o}=\la o', o\ra$. We say:
$\Cn(\bw)$ is the connected component containing the root $o$ and the ends
$\gamma', \gamma$ of the horocyclic product of two percolative subtrees
$T'=T'(\w')$ and $T=T(\w)$ of retention parameter $p'$, and $p$, respectively.\\

{\bf Remarks:} i.)  $T(\w')'\circ T(\w)$ restricted to $\Cn(\bw)$ is the
horocyclic product of two rooted trees (with fixed ends) drawn from
the {\bf augmented Galton Watson measure} [\ref{adalyo}], with
offspring distribution $\{p_k\}$ concentrated on $I:=\{\a_o, ... ,
\a\}$, and being of binomial type 
\bel 
p_k \;=\;\left( {\a - \a_o} \atop k-\a_o \right)p^{k-\a_o}(1-p)^{\a-k+\a_o},
\;\;\;\textrm{ where }k\in I.  \label{eq:law}
\eel

ii.) If $\a_o>0$ or $\a'_o>0$, the percolation $\bar{\mu}$ is {\em not} an
invariant percolation: any vertex-transitive subgroup of the automorphism group
Aut$(G)$ of $G$ must contain the operation which exchanges some vertex
$\bar{k}=\la k', k\ra \in V$ with a vertex $\bar{l}=\la l', l\ra \in V$ at the
same horocycle ($\h(k)=\h(l), \h'(k')=\h'(l')$), where $k$ is connected to its
predecessor (parent) in $\T$ by a marked edge ($\in E_m(\T)$) and $l$ is
connected to its parent by an unmarked edge ($\in E_u(\T)$). Exchanging these
vertices does not leave the measure $\bar{\mu}_{\bare}$ invariant. In
particular, an exchange of these two edges may lead to $\Cn(\bw)$ being
disconnected from one of its fixed ends.  On the other hand, if $\a_o=0,
\a_o'=0$, the model corresponds to an invariant bond-percolation with retention
parameter $p'\cdot p$ (not Bernoulli!). Equivalently, if $\a_o'=\a'$,
$T_{\a_o', \a'}'(\w')$ is deterministic. If under these circumstances $\a_o=0$,
the model is also an invariant percolation, however with retention parameter
$p$. \\ 

\begin{theo}\label{theo:strong_amen}
  Let $p'=p \in[0,1]$. Furthermore, let $\a_o', \a'$, and $\a_o, \a$ be the
  minimum and maximum number of offspring at each site of $T'(\w')=T_{\a'_o,
    \a'}(\w')$ and $T(\w)=T_{\a_o,\a}(\w)$, respectively.  Let $\a_o', \a_o \ge
  1$, and \bel \a_o'\;+\;p'(\a'\;-\;\a_o')\;\;=\;\;\a_o\;+\;p(\a\;-\a_o).
  \label{eq:cond}
\eel

Then the restriction of the horocyclic product
$H(\wb)=T_{\a_o', \a'}\circ T_{\a_o, \a}$ to the connected component
containing the root $\bar{o}$, has $\bar{\mu}$-almost surely an
anchored isoperimetric constant equal to zero, $\bar{\mu}$-almost
surely, i.e.  $H(\wb)|\Cn(\bw)$ is $\bar{\mu}$-a.s. strongly
amenable.\\\label{theo:1}
\end{theo} 

{\bf Remark:} iii.) The strong amenability ($\eqiv$ vanishing anchored
isoperimetric constant [\ref{hss}]) of the Diestel-Leader graphs $DL_{\a,\a}$
is well known: e.g., it follows from the fact that the speed of the simple
random walk is zero [\ref{kv}], together with the conclusion of [\ref{vir}],
that anchored expansion implies positive speed.  $DL_{\a_o',\a_o}$ and
$DL_{\a',\a}$ result as extremal cases $\bar{\mu}$-a.s.  if $p'=p=0$ or
$p'=p=1$, respectively. The theorem says, that amenability is stable under the
`random pertubation' given by the specific construction of the percolation
process, above, if the equal-growth condition (\ref{eq:cond}) is met. \\

{\em Proof:} Let the root $\bar{o}:=\la o', o\ra$ have the level coordinate $0$
(see Fig.1). Let $C_o(\bw)=:T'(\w')\circ T(\w)$, i.e. call $T'(\w')$ and
$T(\w)$ the two random subtrees of $\T'$ and $\T$ with roots $o'$ and $o$,
respectively (see Fig. 2). Let $X_j^{(h)}(\w)$ be the number of leaves (at
level $j$) of the finite subtree $T^{(h)}_j(\w)$ of $T(\w)$, rooted at $-h$ with
height $h+j$. Likewise, call $T'^{(h)}_j(\w')$ the subtree of $T'(\w')$ rooted
at $+h$ with depth $h-j$ ($j\in \{-h, ..., h\}$) and $X'^{(h)}_j(\w')$ its
leaves, also located on level $j$.\\ 

If we find a F\o lner sequence, i.e. a sequence of finite, connected
subgraphs of $T'(\w')\circ T(\w)$ restricted to the connected
component $\Cn(\bw)$ containing the root, with an isoperimetric ratio
converging to zero, then this graph (denoted by $T'(\w')\circ
T(\w)|C_o(\bw)$) is strongly amenable (i.e. has vanishing anchored
isoperimetric constant). In particular, this is true if the finite
subgraphs given by $T'^{(h)}_j(\w')\circ T^{(h)}_j(\w)$ have an
isoperimetric ratio $I_h(\bw)$ as subgraphs of $C_o(\bw)$ , converging
to zero, $\bar{\mu}$-almost surely, as $h\to \infty$. It is clear,
that 
\bel 
I_h(\bw) := \frac{X_h^{,(h)}(\w') +
  X_h^{(h)}(\w)}{\suli{j=-h}{h} X_{-j}^{,(h)}(\w' )
  X_j^{(h)}(\w)}. \label{eq:I_h} 
\eel We show $I_h(\bw)\to 0$,
$\bar{\mu}$-almost surely, to prove the theorem.\\

We note, there is a formular for $X_j^{(h)}$:\\

\begin{lemma} Let $e(k_1, ..., k_l)$ be the edge in the subtree of $\T$ rooted
at $-h$ between the level $l-1$ and $l$, which is uniquely determined by the
$l$-tuple of numbers $k_i\in\{1,2,..., \alpha\}$, with $i\in\{1, ..., l\}$ in
an obvious way: among the $\alpha$ choices of children $k_j$ is chosen on a
path from the root at the $j$th step. Then

\bel
 X_j^{(h)}(\w) = \suli{k_1=1}{\alpha}\cdot\cdot\cdot\suli{k_{h+j}=1}{\alpha}
\prod\limits_{l=1}^{h+j} (\chi_{\{1, 2, ..., \alpha_o\}}(k_l)\;\;+\;\;\chi_{\{
  \alpha_o+1, ..., \alpha\}}(k_l) \chi_{E(\w)}(e(k_1, ..., k_l))).
\label{eq:formular} 
\eel 
\end{lemma}
{\em Proof:} The multiple sum is a sum over the leaves of a homogeneous tree of
finite height, while the product (over $l$) concerns the edges of a path
leading from the root to each of the leaves of these trees. The factors of the
products correspond to indicators of the events of the corresponding
edges being open or closed. \hfill \qed\\

Since $X^{(h)}_j$ is a subtree of a Galton-Watson tree, it is clear
that $\E[X_j^{(h)}]= z^{j+h+1}$, where 
\be 
z:= (1-p)\cdot\alpha_o\;\;+\;\;p\cdot\alpha = \alpha_o +
p(\alpha-\alpha_o).  
\ee 
and $z'$ is the corresponding primed parameter.  This follows also by
using (\ref{eq:formular}) by evaluating the expected value at the
leaves at the highest level, first.\\

In order to simplify the notation used in (\ref{eq:I_h}), let
$X_j'(\w')= X_j^{,(h)}(\w')$ and $X_j(\w)= X_j^{(h)}(\w)$.

We are employing Jensen's inequality in the following way: For any finite sequence
$\{x_j\}_{j=1}^N$,
\bel
\frac{1}{N}\sum_{i=1}^N \frac{1}{x_i} \ge \frac{N}{\sum^N_{i=1} x_i}.\label{eq:jensen}
\eel

If $N=2h+1$, then (\ref{eq:jensen}) applied to (\ref{eq:I_h})  gives

\bel
I_h &\le& \frac{1}{(2h+1)^2} \suli{j=-h}{h} \;
\frac{X'_h + X_h}{X_{-j}' X_j}.
\eel

Define $Y_j:= X_j^{(h)}/ z^{h+1+j}$, and $Y_j':= X_j^{,(h)}/ z^{h-j}$.  As is
well known, $Y_j$ is a martingale [\ref{har}]. Therefore, $\E[Y_{j+1} -
Y_j]=0$, and $\E[X_h/X_j]= z^{h-j}\E[Y_h/Y_j] = z^{h-j}(1\;+
\;\E[(Y_h-Y_j)/Y_j]) = z^{h-j}(1\;+\;\E[Y_h-Y_j]\E[1/Y_j]) =z^{h-j}$.
So, due to the independence between trees, if $\E_\Omega$ is the expectation value
obtained by integration (only) over $\Omega$,

\bel
\E_\Omega\left[ \frac{X_h}{X_{-j}' X_j}\right]&=& \frac{1}{X_{-j}'}
\left(\;z^{h-j}\;\right) = \frac{1}{Y_j'}\,.\label{eq:xy}
\eel

Similarily, denoting by $\E_{\Omega'}$ integration over $\Omega'$, we have
$\E_{\Omega'}\left[ \frac{X'_h}{X_{-j}' X_j}\right] =  \frac{1}{Y_j}$.  It follows,

\bel
\E_\Omega\left[ \frac{X_h}{\sum X_{-j}' X_j} \right] \le
\frac{2}{(2h+1)^2} \suli{j=-h}{h} \frac{1}{Y_j'}, \;\;\;\;\;\;\;\;\;\;\E_{\Omega'}\left[
  \frac{X_h'}{\sum X_{-j}' X_j} \right] \le \frac{2}{(2h+1)^2} \suli{j=-h}{h}
\frac{1}{Y_j},\label{eq:twoces} \eel

for every  $h\in \N$.\\

\eject
\newpage

Since $\E_\Omega Y_j < \infty$, by the martingale convergence theorem,
the sequence $Y_j$ converges $\mu$-almost surey. Moreover, under the
assumption of non-extinction, the probability that $\lim Y_j=0$ is
zero, by the Kesten-Stigum theorem [\ref{hering}]. The same statement
also holds for $\Omega$ and $Y_j$ replaced by $\Omega'$ and
$Y'_j$. This implies

\be 
\wh{C}(\w):=\sup\limits_{j\in \N} \frac{1}{Y_j(\w)} < \infty,\;\;
\mu-\textrm{a.s.}, \;\;\;\;\textrm{ and
}\;\;\;\;\;\wh{C}'(\w'):=\sup\limits_{j\in \N} \frac{1}{Y'_j(\w')} < \infty,
\;\;\;\mu'-\textrm{a.s.}.  
\ee

with $\wh{C}(\w), \wh{C}'(\w')$ independent of $h$.
Using (\ref{eq:twoces}), this implies that

\be 
\E_\Omega\left[ \frac{X_h}{\sum X_{-j}'(\w') X_j} \right] \le
\frac{\wh{C}'(\w')}{(2h+1)}
,\;\;\;\mu'-\textrm{a.s.}
\ee

and

\be
\E_{\Omega'}\left[ \frac{X'_h}{\sum X_{-j}' X_j(\w)} \right]
\le \frac{\wh{C}(\w)}{(2h+1)}, ,\;\;\;\mu-\textrm{a.s.}.  
\ee

Therefore, $\bar{\mu}$-almost surely, 
\be
I_h(\w, \w') = O(h^{-1}),\;\;\;\;\;\;\;\textrm{ as }h\to\infty.
\ee

Since $\bar{\mu}$-almost surely, $\{ I_{h}\}_h$ converges to zero as
$h\to \infty$, the sequence of sub-graphs $T'^{(h)}_j(\w')\circ
T^{(h)}_j(\w)$ restricted to the connected component $\Cn(\bw)$ is a
F\o lner sequence, $\bar{\mu}-a.s.$. Since it is a sequence of
connected sub-graphs each containing the root $\bar{o}= \la o', o\ra$, the
horocyclic product $T'_j(\w')\circ T_j(\w)$ restricted to the
connected component $\Cn(\bw)$ is $\bar{\mu}$-a.s. strongly
amenable. \hfill \qed\\

{\bf Remark:} iv.)  Theorem 2.4,ii. in [\ref{hss}] gives a comparable
result for Bernoulli percolation on transitive graphs.  In our model,
the cases of either $\a_o=0$ and $\a_o'=0$ (both trees' edges all
marked), or $\a_o=0$ and $\a_o'=\a'$ (one tree being deterministic)
imply that all the edges of $DL_{\a',\a}$ are percolative. However,
this is different from Bernoulli percolation in that there are strong
correlations between the openness of the edges $e'\circ e_1$ and 
$e'\circ e_2$.

\eject
\newpage

\section{Amenability of products  of differing percolative 
trees} 

We now proceed to the situation where the `anchor' in the assumptions
is removed and instead of a F\o lner sequence of connected subsets
with a root, only finiteness of the elements of this sequence is
required. Recall that a graph $\tilde{G}$ is amenable iff $\inf
|\partial_{\tilde{G}} C|/|C|=0$, where the infimum is taken over all
finite subsets of $V(\tilde{G})$. As before, $\tilde{G}|C$ will denote
the subgraph of $\tilde{G}$ induced by the subset $C$ of $V(\tilde{G})$.
Once more, let $H(\wb)=\la V, \bar{E}(\wb)\ra$ denote the horocyclic product 
$T'(\w')\circ T(\w)$.\\

For two vertices $u, v\in V(\T)$ we will denote by $u \prec v$ the event
that $v$ is a vertex of a sub-tree of $T(\w)$ with $u$ being the
root. If $u\prec v$ occurrs, this implies that $\h(u) < \h(v)$, and
that $u$ is an ancestor of $v$ in the hierarchy induced by the
Busemann-function $\h$. Similarily, for $u', v'\in V(\T')$, we define
$u'\prec v'$ to mean $\h'(u') < \h'(v')$, and $u'$ is an $'v'$ ancestor
in the hierarchy induced by $\h'$.\\

\begin{theo}\label{theo:2}
  Let \,$H(\wb)=T'(\w')\circ T(\w)$ be a horocyclic product of two
  trees with offspring-measures which have non-empty support.,
  i.e. let $\{ \a'_o, ..., \a\} \cap \{\a_o, ..., \a\} \neq \emptyset$.
  Let $p', p \in (0,1)$.  Then $H(\wb)|\Cn(\bw)$ is almost surely amenable.
\end{theo}

Consider at first the horocyclic product of two finite
percolation sub-trees with equal number of minimal offspring $\a'_o=\a_o$,
maximal offspring $\a'=\a$, and equal height $N\in\N$.
Choose $\bw\in\Omega'\times\Omega$ such that
\bel 
\textrm{deg}(\bar{v}) = 2\alpha_o.\label{eq:an} 
\eel

This is the event that all percolative edges are closed. By definition
of $\bar{\mu}$ (see (\ref{eq:law})), the probability that this occurrs
on both trees is $(1-p')^{2M_N}(1-p)^{2M_N}>0$, where $M_N=
\a\cdot(\alpha_o^{N+1}-1)/(\a_o-1)$, the factor after $\a$ being the
number of vertices in a $\a_o$-regular tree of height $N$.  Likewise,
the event of each vertex in either of two rooted sub-trees within a certain
finite intervall of levels having the same number $\beta$ of offspring
has also positive probability. Due to the invariance of $\bar{\mu}$
under shifts accross different horocycles (=levels), any arbitrary
number of such events occurrs also with positive
probability. Moreover, due to the independence of the openness of
edges between vertices on different horocycles, the measure is ergodic
with respect to this shift. \\

{\em Proof:} (Theorem 3.1) Since $\a_o\ge 1$, the tree $T(\w)$ is
infinite. Choosing any sequence $(v_j)_{j\in \N}$ of vertices in
$V(T(\w))$ with $v_{0}=o$, and ${\bf h}(v_n)=n$, we may consider the
induced sequence of vertices $(\bar{v}_n)_{n\in \N}$ with
$\bar{v}_n=\la v_n', v_{n}\ra \in \Cn$, and $v'_0=o' \in
V(T'(\w'))$. Note that $\h'(v'_n)=-n$, and that therefore
$v'_{j+N}\prec v'_n$, for every positive $N\in\N$. \\

Let $Z_v(\w)$ be the number of offspring of $v\in V(T(\w))$, and
$Z'_{v'}(\w')$ the number of offspring of $v'\in V(T'(\w'))$. For
given $\beta \in \{ \a'_o, ...., \a'\}\cap \{ \a_o, ...., \a\}$, $N\in
\N$, and $j\in \N$, the event
\be
C^N_j := \left\{ \la \w', \w\ra \in
\bar{\Omega}\;:\; Z'_{u'}(\w') = Z_u(\w) = \beta
\;\textrm{for } \la u', u\ra \in \Cn(\wb) \; \textrm{ s.t. }
\; v_j\prec u, \;\;v'_{j+N}\prec u'  \right\} 
\ee

occurrs with positive probability $\bar{\mu}(C^N_{j})>0$.\\

The measure $\bar{\mu}$ is ergodic under the shift $\bar{v_k}\mapsto
\bar{v}_{k+1}$ because of the independence of the openness of
different edges (in particular of edges on different horocycles).
Therefore, there is an $n\in \N$ for each $N\in\N$, such that
$C^N_n$ occurrs. Defining the random sequence $n_\cdot(\bar{w})\in \N^\N$ by
\be
n_N(\wb) := \inf\left\{ k\in N\; :\; \bw\in C^N_{k}\;\right\},
\ee

we identify a sequence of vertex sub-sets (called {\em tetraeder} in
[\ref{barwoe}]):
\be
V_N\{ \la u', u\ra \in 
\Cn(\wb)\;:\; v_{n(N)} \prec u, \textrm{ and }v'_{n(N)+N}\prec u'\}.
\ee 
This sequence is a F\o lner-sequence, since the isoperimetric constant of the
finite subgraph of $H(\bar{w})$ induced by $V_N(\wb)$ is given by

\be
I_H(H(\bw)|V_N(\bw)) = \frac{|\partial_{H(\bw)}V_N(\bw)|}{|V_N(\wb)|} 
= \frac{2\beta^N}{\sum_{j=0}^N \beta^j \beta^{N-j}}=\frac{2}{N+1}.
\ee

 \hfill\qed\\

{\bf Remarks:} vi.) In this proof, for the construction of the F\o lner
sequence it is important that with positive probability the graph
locally contains arbitrarily large but finite subgraphs of the graphs
induced by $C_o(\bw)$, which are (finite) symmetric horocyclic products.\\

 vii.) The proof can be transferred to the situation in which
 dependent percolation prevails, however, with independence between
 the different trees and stationarity and stationary ergodicity with
 respect to the shift between the different horocycles (levels) within
 a single tree.\\

\eject
\newpage

\section{Non-amenability}\label{sec:purenona}

One difficulty in bounding the anchored isoperimetric constant of subgraphs of
horocyclic products is, that these subgraphs need not be horocyclic products,
themselves. We will overcome this difficulty by removing additional edges from
$H(\bw)$ and recognising the remaining graph to have a uniformly bounded
isoperimetric ratio in a deterministic, non-amenable horocyclic product
(compare with [\ref{rau}], [\ref{mare}]). Theorem 1.1 in [\ref{cheper}] and
Theorem 2.4,i. in [\ref{hss}] refer to Bernoulli percolation on a not
necessarily transitive but locally finite graph. Again, this situation is
realised in our model if $\a_o=0$ and $\a_o' \in \{0, \a'\}$.\\

Using the idea, that the removal of edges may lead to non-amenable subgraphs,
we formulate the following lemma. (These are two statements one with, the other
without the parentheses.)

\begin{lemma}\label{lemma:nona}
  If, given a graph $G=\la V, E\ra$ of bounded degree, there is a subset of
  edges $E'\subset E$ , such that for every finite, connected induced subgraph
  $G|C$ with $C\subset V$ (containing the root), the connected components of
  $G\sme E'\;\; |\; C$ each have an isoperimetric ratio in $G\sme E'$ uniformly
  bounded from below by $i_o>0$, $G$ is (weakly) non-amenable with (anchored)
  isoperimetric constant greater or equal to $i_o$.
\end{lemma}

{\em Proof:} By lemma A3.3 of [\ref{hss}], since the graph $G$ has bounded
degree, it is sufficient to consider only {\em connected} subgraphs $G|C$ in
the assumption (for non-amenability - for weak non-amenability, 'connected'
requires no further justification).  Let $\{C_j\}$ with $C_j \subset C$ be the
finite set of disjoint subsets for which $(G\sme E') | C \;=\;\cup_j (G\sme E')
| C_j=\sum_j G|C_j$. In other words, after taking away the edges $E'$ of $G|C$,
we are left with the disjoint subgraphs $G|C_j$.  Since by assumption, each of
them fulfills $|\partial_{G\sme E'}C_j| \ge i_o|C_j|$ for some positive $i_o$
(independent of $C$), it holds that
\be
\frac{|\partial_G C|}{|C|} \ge \frac{|\partial_G C \setminus E'|}{|C|} =
\frac{|\partial_{G\sme E'} C |}{|C|} = \frac{1}{|C|}\sum_j
|\partial_{G\sme E'} C_j| \ge \frac{i_o}{|C|} \sum_j |C_j| = i_o. \hfill \qed
\ee

What if the smallest possible number of offspring of one of the trees ($\a_o$)
is larger than the largest number of offspring of the other tree ($\a'$)?  We
answer this question in the case of one tree being deterministic $(\a_o'=\a'$).

\begin{theo}
  Let $\a_o'=\a'<\a_o<\a$. Then there is non-amenability, {\em for all}
  realisations of the random subgraphs $\T'_{\a'}\circ T_{_{\a_o, \a}}(\w)$ of
  DL$_{\a', \a}$.
\end{theo}

{\em Proof:} Note that when all the percolative edges $E_p$ are removed
$(p=0)$, the remaining graph is disconnected and all its connected components
infinite and non-amenable. The isoperimetric ratio of any connected sub-graph
induced by a finite subset of verticies $W\subset V$ fulfills 
\be
\frac{|\partial_H W|}{|W|} \ge \frac{|\partial_{H\sme E_p} W|}{|W|}, 
\ee

which is of the form $(\suli{j}{}|\partial_{H\sme E_p} W_j|)/ \suli{j}{}|W_j|$,
for a finite number $n$ of subsets $W_j$ of $W$, where the graph induced by
each $W_j$ is, by Lemma \ref{lemma:first}, the connected component resulting
from removing the percolative edges $E_p$.  Note that each of the corresponding
isoperimetric ratios $I_j=|\partial_{H\sme E_p} W_j|/|W_j|$ is uniformly
bounded away from zero, since it is a subgraph of DL$_{\a', \a_o}$, which is
non-amenable, and by assumption $\a'<\a_o$.  Given $n$ ratios $|\partial
W_j|/|W_j|\;\ge\;c$, uniformly bounded by $c>0$, we have that $\suli{j}{}
|\partial W_j| \;/ \suli{j}{} |W_j| \ge c.$ This means that {\em every}
realisation of $H$ is
non-amenable. \hfill \qed\\

{\bf Remark:} viii.) Note that a finite graph has vanishing isoperimetric
constant.

\section{Summary, Outlook and 
Acknowledgements}

In this paper, it was proven that there are transitive graphs with certain
independent percolation processes for which either strong amenability or 
amenability prevails, depending on the choice of the retention parameters.
Two methods have been introduced to investigate amenability of percolative
subgraphs: 1.) The expected isoperimetric ratio may lead to the existence of
a F\o lner sequence, and 2.) the removal of edges may allow comparison with
random graphs for which non-amenability has been proven.\\

{\bf Question:} Is the range of strong amenability restricted to the assumption
given in 

Theorem  \ref{theo:strong_amen}, or is there a non-trivial regime of the
parameters $p', p$ for anchored expansion (weak non-amenability) to prevail?\\

The question wether for some horocyclic products of trees drawn from the
augmented Galton-Watson measure there is (strong) amenability together with
simple random walk having positive speed is answered in a forthcoming paper of
the author with
V. Kaimanovich.\\

The question to what extent similar results hold for non-random periodic and
quasi-periodic trees is investigated in a project by the author, D. Lenz, and
I. Veseli\`c.\\

I am grateful for helpful discussions with V.Kaimanovich, S. Brofferio., W.
Woess. I owe special thanks to T.  Antunovic, S. M\"uller, Rainer
Siegmund-Schultze, and E. Candellero for hints and corrections, and to
P. Mathieu who observed theorem \ref{theo:2}. I also thank Y. Zhang for
pointing out the role of the finite subtrees in an earlier version of this
paper. His remark lead to the discovery of a mistake in an attempt to prove
weak non-amenability. A characterisation of anchored expansion of horocyclic
products is given for a class of trees with stronger growth-conditions in a
paper, which is joint work with D.Lenz and I. Veseli\'c.

\section{Bibliography}

\begin{enumerate}

\item{S.Adams, R.Lyons: `Amenability, Kazhdan's property and percolation for
    trees, groups and equivalence relations', Israel Journ., Math. 75, 341-370,
    1991}\label{adalyo}

\item{L.Bartholdi, W.Woess:`Spectral computations on lamplighter groups and
    Diestel-Leader graphs',  J. Fourier Analysis Appl., 11, 2, 175--202, 2005
  }\label{barwoe} 

\item{I. Benjamini, R. Lyons, O. Schramm: `Percolation Pertubations in
    Potential Theory and Random Walks', }\label{bls}

\item{D.I.Cartwright, V.A.Kaimanovich, W.Woess: `Random walks on the affine
    group of local fields and of homogeneous trees', Ann. de l'inst. Fourier,
  tome 44, $n^o$ 4(1994), p. 1243-1288}\label{cakawo}

\item{D. Chen, Y. Peres, G. Pete: `Anchored expansion, percolation and speed',
    Ann. of Probab. 32, {\bf no.} 4, 2978-2995, 2004}\label{cheper} 

\item{D. Chen, Y. Peres: `The Speed of Simple Random Walk and Anchored
    Expansion on Percolation Clusters: an Overview', Discrete random walks
    (Paris, 2003), Discrete Math. Theor. Comput. Sci. Proc., AC, 39-44, 2003
    }\label{cheper1} 

\item{R. Diestel, I. Leader: `A conjecture concerning a limit of non-Cayley
    graphs', J.Algebraic Combin. 14(2001), no.1, 17--25}\label{dl}

\item{J. Dodziuk: `Difference equations, isoperimetric inequality and
    transience of certain random walks', Trans. Amer. Math. Soc.  284, no. 2,
    787--794, 1984}\label{dod}

\item{A. Eskin, D. Fisher, K. Whyte: `Coarse differentiation of quasi-isometries
    I: spaces not quasi-isometric to Cayley graphs',   Preprint }\label{efw}

\item{T.E. Harris: `The theory of Branching Processes', Springer, chap. I,
    Theor. 8.3, 1963 }\label{har}

\item{O. H\"aggstrom, R.H. Schonmann, J.  Steif: `The Ising model on
    diluted graphs and strong amenability', Ann.  Probab. 28, no. 3,
    1111--1137, 2000}\label{hss}

\item{F. den Hollander: `Large Deviations', Fields Institute Monographs, AMS, 
chapter I.2, 2000 }\label{denhol}

\item{ Kaimanovich, V. A.; Vershik, A. M.: `Random walks on discrete groups:
    boundary and entropy',  Ann. Probab.  11  (1983),  no. 3,
    457--490}\label{kv}

\item{V. A. Kaimanovich: `Equivalence relations with amenable leaves
need not be amenable', Topology, Ergodic Theory, Real Algebraic Geometry,
Amer. Math. Soc. Transl. (Ser. 2) 202, 151-166, 2001}\label{vadim} 

\item{V. Kaimanovich, F.Sobieczky: `Random walks on horocyclic products',
    preprint, 8, 2007}\label{vaflo}

\item{Asmussen, S.; Hering, H.: 
`Branching processes',
 Progress in Probability and Statistics, 3. Birkh\"auser Boston, Inc., Boston,
 MA, 1983 }\label{hering} 

\item{R. Lyons: `Phase transitions on nonamenable graphs', J. Math. Phys. 41,
    1099--1126, 2000 }\label{lyo}  

\item{R. Lyons, R. Pemantle, Y. Peres: `Ergodic theory on Galton-Watson trees:
    speed of random walk and dimension of harmonic measure', Erg. Theory
    Dyn. Systems 15 , 593-619, 1995}\label{lpp}

\item{P. Mathieu, E. Remy: `Isoperimetry and heat kernel decay on percolation
    clusters', Annals of Prob., 32 1A, 100-128}\label{mare}

\item{B. Mohar: `Isoperimetric inequalities, growth, and the spectrum of
    graphs', Lin. Alg. Appl. 103, 273-285, 1988}\label{mohar} 

\item{C. Rau: `Sur le nombre de points visite\'es par une marche al\'eatoire
    sur un amas infini de percolation´, (preprint), arXive:math.PR/0605056,
    2006}\label{rau}

\item{M. Soardi, W. Woess: `Amenability, unimodularity, and the
    spectral radius of random walks on infinite graphs',
    Math. Zeitsch. 205, 471-486, 1990}\label{sowo}

\item{C. Thomassen: `Isoperimetric inequalities and transient random walks on
    graphs', Ann. of Probab. 20, No.3 1592-1600, 1992 }\label{tho}

\item{Virag, B.:`Anchored expansion and random walk', Geom. Func. Anal. 10,
    1588-1605, 2000}\label{vir}

\item{D. West: `Introduction to Graph theory', Prentice-Hall, 1996, (Definition
    1.2.14)}\label{wes}

\item{W. Woess: `Random Walks on Infinite Graphs and Groups', 
    Cambridge Tracts in Mathematics 138, Camb. Univ. Press, 2000,
    (chapter 12.13)}\label{woess} 

\item{W.Woess: `Lamplighters, Diestel-Leader graphs, random walks, and harmonic
    functions', Combinatorics, Probability \& Computing 14, 415-433, 2005
  }\label{woess2}

\end{enumerate}

\end{document}